\documentstyle[11pt]{article}

\title{\Large \bf  A note on Closure Operators in Category of Modules}

\author{Vishvajit V S Gautam\thanks{vishvajit@imsc.res.in, gautamvvs@yahoo.com }\\
\sl  The Institute of Mathematical Sciences\\
\sl  CIT Campus, Taramani, Chennai - 600113 INDIA}
\date{ }
\begin{document}
\newtheorem{defi}{\bf Definition}[section]
\newtheorem{th}[defi]{\bf Theorem}
\newtheorem{pro}[defi]{\bf Proposition}
\newtheorem{lemma}[defi]{\bf Lemma}
\newtheorem{coro}[defi]{\bf Corollary}
\newtheorem{ram}[defi]{\bf Remark}
\newtheorem{pn}[subsubsection]{\bf Proposition}
\newtheorem{zm}[subsubsection]{\bf Theorem}
\newtheorem{cy}[subsubsection]{\bf Corollary}
\newtheorem{df}[subsubsection]{\bf Definition}
\newtheorem{rk}[subsubsection]{\bf Remark}
\maketitle
\begin{abstract}
In this article  we give application of closure operators in category 
of modules. Our main result shows that every subcategory $\cal A$ of
injective modules of {\bf R-mod} (under a mild condition) 
induces a torsion theory of {\bf R-mod}.\\


\end{abstract}

Notion of  closure operators(operations,  systems,  functions, 
relations) is known to us from algebra. logic, lattice theory and topology.
Categorical view of closure operators play an important role in various
branches of mathematics. In an arbitrary category $\cal X$ with suitable
axiomatically defined notion of subobjects, a (categorical) 
closure operator $c$ is defined to be a family $(c_{X})_{X \in {\cal X}}$ 
satisfying the properties of extension, monotonicity and continuity.
Concept of closure operators induced by a subcategory $\cal A$ of {\bf Top} was 
introduced by Salbany in 1976 and have been studied by many researchers in a 
variety  of situations. Closure operators are proved to be
useful in study of Galois equivalence between certain factorization systems. 
In category of {\bf R-mod} of R-modules closure operators correspond to 
preradicals. For more details, see \cite{[2]}.

In section 2 of this paper, we give characterizations of {\bf R-mod} homomorphisms in terms of
regular  closure operators. Using these characterizations we establish
that every subcategory $\cal A$ of injective modules of {\bf R-mod} (under a mild condition) 
induces a torsion theory of {\bf R-mod} [Theorem 2.3]. Some useful 
observations of  these characterizations to modules over a {\em free ideal rings}({\em semifir or n-fir}) are given in section 3.\\
Subcategories are always assumed to be full and isomorphism closed.

\section{Preliminaries}
Throughout this paper we consider a category $\cal X$ and a fixed class $\cal M$ 
of monomorphisms in $\cal X$ which contains all isomorphisms of $\cal X$. It is 
assumed that\\
$\bullet$ $\cal M$ is closed under composition;\\
$\bullet$ $\cal X$ is finite $\cal M$-complete.\\

A {\em closure operator} $c$ on the category $\cal X$ with respect to class 
$\cal M$ of subobjects is given by a family $c = (c_{X})_{X \in \cal X}$ of 
maps $c_{X} : {\cal M}/X \longrightarrow {\cal M}/X$ such that for every 
$X \in \cal X$ \\
1.  $m  \leq  c(m)$;  2.   $m \leq m' \Rightarrow  c (m) \leq c (m')$;
and  3. for  every $f:X \longrightarrow Y$ and  $m \in {\cal M}/X$,
\hspace{.5in} $ f(c_{X}(m)) \leq c_{Y}(f(m))$.\\
\\
For each $m \in {\cal M}$  we denote by $c(m)$ the c-closure of $m$.\\
An $\cal M$-morphism   $m \in {\cal M}/X$ is  called $c$-closed
if  $m \cong c_{X}(m)$. A closure operator $c$ is said to be idempotent
if  $c(c(m)) \cong c(m)$. In case $ c(m \vee n) \cong c(m) \vee c(n)$ we say 
$c$ to be additive. An $\cal M$-subobject $m$ of $X$ is called $c$-dense in
$X$ if $c_{X}(m) \cong 1_{X}$.
\\
For a subcategory $\cal A$ of $\cal X$, a morphism $f:X \longrightarrow Y$
is an $\cal A$-regular monomorphism  if it is the equalizer of two morphisms
 $h,k:Y \longrightarrow A$  with  $A \in {\cal A}$.

 Let   $\cal M$  contain  the  class  of  regular monomorphisms of $\cal X$.
 For $m:M \longrightarrow X$ in $\cal M$ define 
 \[c_{\cal A}(m) = \bigwedge \{ r \in {\cal M} \mid r \geq m \,\,\,\, and\,\,
 r \,\,\,is \,\, {\cal A}-regular \} \]
which  is  a closure  operator  of $ {\cal X}$.  These closure operators
are called regular  and $c_{\cal A}(m)$  is  called  the  $\cal A$-closure
of $m$. In case ${\cal A} = {\cal X}$ we denote $c_{\cal A}(m)$ by $c(m)$.

\section{Closure Operators in category of modules} 
In this section we will see application of closure operators in category of
modules.\\
Let ${\cal X} =$ {\bf R-mod} the category of $R$-modules and let $\cal M$ be the
class of all monomorphisms of $\cal X$. In this Case $\cal X$ is 
$\cal M$-complete.\\
Recall that by a linear function $f$ on an $R$-module $M$ we mean an R-module
homomorphism $f:M \longrightarrow R$.\\
Let $\cal A$ be a subcategory of {\bf R-mod}. For an object $M$ of {\bf R-mod}, 
${\cal M}/M$ is the set of all submodules of
$M$, we can identify  each $n$ in ${\cal M}/M$ with $N \subseteq M$ a submodule
of $M$. Let $N$ be a submodule of $M$. The $\cal A$-regular closure operator 
of $N$ can be computed by following formula (\cite{[2]})
\[ c_{\cal A}(N) = \cap \lbrace  Ker(g) \mid g\,\,\,  is\,\,\, an\,\,\, 
R-module \,\,\, homomorphism \] \[ from\,\,\ M\,\,\, to\,\,\, 
A \in {\cal A} \,\,\, and\,\,\,g(N) = (0)\rbrace .\]

{\em Throughout the remainder of this paper  $\cal A$ denotes a category of injective
 modules}.
\begin{pro}
 $c_{\cal A}(N) = M$ i.e. $N$ is $c_{\cal A}$-dense in $M$
if and only if\\
 $Hom_{R}(M/N, A) = 0$, $A \in \cal A$.
\end{pro}
{\bf Proof.} Let $A$ be an object in $\cal A$. Consider the exact sequence
\[0 \longrightarrow N \longrightarrow M \longrightarrow M/N \longrightarrow 0\]
Since $Hom_{R}( -, A) $ is left exact functor, the sequence 
\[0 \longrightarrow Hom_{R}(M/N, A) \longrightarrow Hom_{R}(M, A)
\longrightarrow Hom_{R}(N, A)\longrightarrow 0 \] is exact. If $Hom_{R}(M/N, A) = 0$, above
left exact sequence implies 
$$Hom_{R}(N, A) = Hom_{R}(M, A).$$ 
Therefore every $f \in Hom_{R}(M, A)$
satisfying $f(N) = 0$ must be zero. Hence $c_{\cal A}(N) = Ker(0) = M$.\\
Conversely, assume that $f \in Hom_{R}(M/N, A)$ for some $A \in \cal A$.
Let $\pi : M \longrightarrow {M/N}$ be the canonical projection. The
composition $ f \cdot \pi : M \longrightarrow A$ is an R-module homomorphism
such that $f(\pi (n)) = f(N) = 0$ for all $n \in N$. This gives
$c_{\cal A}(N) \subseteq Ker(f \cdot \pi )$ and $M = Ker(f \cdot \pi)$,
which implies $f \cdot \pi = 0$, and so $f(\pi (M)) = 0(M) = 0$ implies
$f(M/N) = 0$ and hence $f = 0$. $\Box$

\begin{pro}
 $c_{\cal A}(N) = N$  i.e. $N$ is $c_{\cal A}$-closed in $M$ if and only if \\
$Hom_{R}(T/N, A) \not= 0$ for each non-zero submodule $T/N$ of $M/N$ and
$A \in \cal A$.
\end{pro}
{\bf Proof.} Let $A \in \cal A$. If $c_{\cal A}(N) \not= N$, we can regard
${c_{\cal A}(N)}/N$ as a non-zero submodule of $M/N$. Suppose 
$0 \not= g \in Hom_{R}({c_{\cal A}(N)}/N, A)$. Using the arguments similar to
the previous Proposition we get $g({c_{\cal A}(N)}/N) = 0$, which implies
$g = 0$. A contradiction, therefore we must have $c_{\cal A}(N) = N$. \\
Conversely, if $c_{\cal A}(N) = N$, $c_{\cal A}(N)$ is not dense in $M$.
Therefore by above proposition we have $Hom_{R}(M/N, A) \not= 0$. Let
$T/N$ be a non-zero submodule of $M/N$. 
Let $0 \not= f \in Hom_{R}(M/N, A)$. Then restriction of $f$ to $T/N$ gives
$Hom_{R}(T/N, A) \not= 0$. $\Box$\\
\\
Next we relate the $\cal A$-regular closure operator defined in the beginning
of this section with torsion theory of {\bf R-mod} in sense of \cite{[3]}.

\begin{th}
Subcategory $\cal A$ of {\bf R-mod} satisfying ${\cal T} \cap {\cal A} = $ {0} for a class ${\cal T}$ of objects 
of {\bf R-mod} induces a torsion theory 
for {\bf R-mod} such that $\cal A$ forms a torsion free class consisting
of torsion free objects.
\end{th}
{\bf Proof.} Consider the pair $({\cal T}, {\cal F})$ of classes of
suboject of {\bf R-mod}, where
$$ {\cal T} = \lbrace M/N \mid c_{\cal A}(N) = M \rbrace $$ and
$$ {\cal F} = \cal A $$
It is immediate from above two propositions that $({\cal T}, {\cal F})$
is a torsion theory.\\
This torsion theory is also a {\em hereditary torsion theory}. $\Box$

Let $({\cal T}, {\cal F} = {\cal A})$ be a torsion theory induced by
a subcategory $\cal A$ of {\bf R-mod} in sense of Theorem 2.3 . Following
result (cf. [3], Proposition 2.1, 2.2) characterizes the 
torsion classes.
\begin{th}
For torsion theory $({\cal T}, {\cal F} = {\cal A})$
\begin{enumerate}
\item The class $\cal T$ of objects of {\bf R-mod} is closed under image, infinite
sums and group extensions.
\item $\cal A$ is closed under kernel, infinite products and group 
extensions.
\end{enumerate}
\end{th}

\section{Closure Operators in category of firs ( semifirs or $n$-firs)}

Recall that by a {\em free ideal ring} we mean a ring $R$ with the property
that all their (one-sided) ideal are free as left module (resp. right module)
of unique rank. \\

Let {\bf R}-{\em mod} denotes the category of modules whose underlying
rings are {\em firs (semifirs or n-firs)} (\cite{[1]}).  Let $R$ be injective as $R$-module.\\
An $R$-module $M$ is said to be bounded if $Hom_{R}(M, R) = 0$
\begin{pro}
Let $N$ be a submodule of $M$.
\begin{enumerate}
\item A quotient module $M/N$ is bounded if and only if $c_{\cal A}(N) = M$.
\item A quotient module $M/N$ is unbounded if and only if $c_{\cal A}(N) = N$.
\end{enumerate}
\end{pro}
{\bf Proof.} Obvious. $\Box$\\
Clearly bounded injective modules form a torsion class in above sense.
\begin{pro}
 
\begin{enumerate}
\item $c_{\cal A}(N) = N$ if and only if $M/N$ has a direct
summand which is free of positive rank.
\item $c_{\cal A}(N) = M$ if and only if $M/N$ does not have a direct
summand which is free of positive rank.
\end{enumerate}
\end{pro}
{\bf Proof.}\,\,1. \, Since the image of every non-zero homomorphism $f$ is
free, which implies $M/N$ has a direct summand.\\
Conversely, If $M/N$ is a direct summand of a free module $f(M/N)$ of 
positive rank. Then $f$ should be free.\\
2. \,\, The result is immediate from 1. $\Box$\\
\\
Many fragment of torsion theory for {\bf R-mod} can be related with
results proved in this section. Theorem 2.3  provides a convenient way to
do this. One can use these result in the study of Gabriel topology \cite{[3]}. 
If we assume R to be {\em fir}({\em semifir}, {\em n-fir} or Bezout domain),
there are some areas of {\em free ideal ring theory} where these results
are useful. It may be a subject of further investigation.\\
\\
\\

\end{document}